%% file: central2.tex
\documentclass{article}   
\usepackage{amsmath}   
\usepackage{amsthm}
\usepackage{epsfig}
\usepackage{amssymb} \textwidth=6in \hoffset=-1.3cm
\theoremstyle{plain}
\newtheorem{theorem}{Theorem}[section]
\newtheorem{lemma}[theorem]{Lemma}
\newtheorem{proposition}[theorem]{Proposition}
\newtheorem{corollary}[theorem]{Corollary}

\theoremstyle{definition}
\newtheorem{definition}[theorem]{Definition}
\newtheorem{example}[theorem]{Example}

\theoremstyle{remark}
\newtheorem{remark}[theorem]{Remark}

\title{Computation of Centralizers in Braid groups and Garside Groups\footnote{Both authors
partially supported by the European Network TMR Sing. Eq. Diff. et Feuill.}}
\author{Nuno Franco\footnote{Partially supported by SFRH/BD/2852/2000.} \and
Juan Gonz\'alez-Meneses\footnote{Partially supported by MCYT, BFM2001-3207 and FEDER.}}
\date{December, 2002}
\begin{document}
\maketitle

\begin{abstract}
  We give a new method to compute the centralizer of an element in Artin braid groups
and, more generally, in Garside groups. This method, together with
the solution of the conjugacy problem given by the authors in
\cite{FGM}, are two main steps for solving conjugacy systems, thus
breaking recently discovered cryptosystems based in braid
groups~\cite{AAG}. We also present the result of our computations, where we notice that
our algorithm yields surprisingly small generating sets for the centralizers.
\end{abstract}

\vspace{.5cm}
\noindent This paper is dedicated to Jos\'e Luis Vicente C\'ordoba, on his $60^{\mbox{th}}$ birthday.

\section*{Introduction}
Given a group $G$, the {\em centralizer} of an element $a\in G$,
denoted $Z(a)$, is the subgroup of $G$ consisting of all elements
which commute with $a$. Our goal in this paper is to give a good
algorithm to compute a generating set for the centralizer of an
element in a {\em Garside group}.

Garside groups were introduced by Dehornoy and Paris~\cite{DP} (their original name was
{\em small Gaussian groups}, but there has been a convention to call them Garside groups).
We will consider Artin braid groups~\cite{A} as the main examples of Garside groups. Given an integer
$n\geq 2$, the {\em braid group on $n$ strands}, $B_n$, is defined by the following presentation:
\begin{equation}
B_{n}=\left\langle
\begin{array}{cc}
\sigma _{1},\sigma _{2},\ldots ,\sigma _{n-1} & \left|
\begin{array}{ll}
\sigma _{i}\sigma _{j}=\sigma _{j}\sigma _{i} & (|i-j|\geq 2) \\
\sigma _{i}\sigma _{i+1}\sigma _{i}=\sigma _{i+1}\sigma _{i}\sigma _{i+1} &
(i=1,\ldots ,n-2)
\end{array}
\right.
\end{array}
\right\rangle .  \label{presen}
\end{equation}

Braid groups are of interest not only in Combinatorial Group Theory, but also in Low Dimensional Topology and,
more recently, in Cryptography. Other examples of Garside groups are spherical (finite type) Artin groups \cite{BS}
and torus knot groups, among others.

Computing centralizers in a Garside group is of interest in
itself, but can also be applied to solve other questions. For
instance, consider two elements $a,b$ in a Garside group $G$.
Suppose that we know an element $c\in G$ that conjugates $a$ to
$b$, that is, $c^{-1}ac=b$. Consider then the set
$Z_{a,b}=c\:Z(b)=\{c\alpha:\; \alpha\in Z(b)\}\subset G$. Then
$Z_{a,b}$ is the set of all elements in $G$ that conjugate $a$ to
$b$: Indeed, an element $d\in G$ conjugates $a$ to $b$ if and only
if $d^{-1}ad=b$, then
$b=d^{-1}(cc^{-1})a(cc^{-1})d=(d^{-1}c)b(c^{-1}d)$, so $c^{-1}d\in
Z(b)$; hence $d\in Z_{a,b}$.

 This property may be used for solving {\em conjugacy systems} in Garside groups:
Given $a_1,a_2,\ldots, a_k,b_1,b_2,\ldots,b_k\in G$, find an element $c\in G$ such that
$c^{-1}a_ic=b_i$, for $i=1,\ldots,k$. The solutions of such a system are the elements in
$Z_{a_1,b_1}\cap \cdots \cap Z_{a_k,b_k}$. These kind of problems play a central role in some new public-key
cryptosystems (see \cite{AAG} and \cite{K-L}), based on braid groups. To break such cryptosystems, one must solve a
conjugacy system such as the previous one.

 The conjugacy problem in braid groups has been solved by Garside \cite{G}, and his algorithm has been improved in~\cite{M} and generalized
to all Garside groups in~\cite{Pi2}. In~\cite{FGM}, the authors gave a more efficient algorithm than all the above, to
solve the conjugacy problem in all Garside groups. So, given two conjugated elements $a,b\in B_n$, we know how to find an
element $c\in G$ such that $c^{-1}ac=b$. Using the
algorithm that we shall explain in this paper, we can compute a
generating set for $Z(b)$, hence we know how to generate elements in
$Z_{a,b}$. We still do not know how to compute an element in
$Z_{a_1,b_1}\cap \cdots \cap Z_{a_k,b_k}$ even if we know how to
generate elements in each $Z_{a_i,b_i}$. We believe that a deeper
study of the structure of centralizers in Garside groups will
provide a solution to this problem. Anyway, we think that the
algorithm we give to compute centralizers in Garside groups is a
good step towards the solution of these systems.

There exists another algorithm to compute the centralizer of an
element in braid groups, which was given by Makanin~\cite{Ma}. It
can be easily generalized to all Garside groups, but it is a
fairly theoretical algorithm, which has a huge complexity and
gives a large amount of redundant generators. One
could also make use of the bi-automatic structure of Garside groups~\cite{D}
to find the centralizer of an element. But this also seems quite inefficient.
The new method that we introduce is quite simple and surprisingly efficient. Actually,
the generating sets obtained in our computations with braid groups
are so small, that they led us to conjecture that the
centralizer of any braid in $B_n$ can be generated by no more than $n-1$ elements.

 After writing an early version of this paper, we were told by M. Korkmaz of a family of counterexamples
to this conjecture, due to N. V. Ivanov (the smallest counterexample belongs to $B_9$, while our
computations were up to $B_8$). Nevertheless, it has been recently proven by the second author and Bert
Wiest~\cite{GW} that, for $a\in B_n$, $Z(a)$ can be generated by less than $\frac{n(n-1)}{2}$ elements.

  The algorithm in this paper works as follows: given an element $a$ in a Garside group $G$, it constructs a graph $\Gamma$ associated
to $a$, such that the fundamental group of $\Gamma$ maps onto
$Z(a)$. Then it computes a generating set for the fundamental
group of $\Gamma$, which maps to a generating set for $Z(a)$.

This paper is structured in the following way: In Section~1 we give the basic definitions and results concerning Garside groups.
In Section~2, we introduce a special kind of elements, the {\em minimal simple elements}, which are used to construct the graph
$\Gamma$. This graph is studied in Section~3. We explain our algorithm in detail in Section~4, then we study its complexity in
Section~5 and, finally, in Section~6 we present the results obtained by implementing the algorithm.

\section{Garside groups and simple elements}\label{secsimple}

In this section we will give the definitions of Garside monoids and groups, and the basic results
which we shall need to present our algorithm. To find the proofs of the results, and more details,
see \cite{G}, \cite{M}, \cite{T}, \cite{BKL}, \cite{DP}, \cite{D} and \cite{Pi}.

 Consider a cancellative monoid $M$, with no invertible elements. We can define a partial order on its elements, called the
{\em prefix order}, as follows: For $a,b\in M$, we say that $a\prec b$ if $b$ can be written in such a way that $a$ is a prefix
of $b$, that is, if there exists $c\in M$ such that $ac=b$. In this case, we say that $a$ is a left
divisor of $b$. There also exists the {\em suffix} order, but we will not use it in this paper, so in the
above situation we will just say that $a$ divides $b$, or that $b$ is a multiple of $a$.

 Given $a,b\in M$, we can naturally define their (left) {\em least common multiple}, $a\vee b$,
and their (left) {\em greatest common divisor}, $a\wedge b$, if they exist.
That is, $a\vee b$ is the minimal element (with respect to $\prec$) such that
$a\prec a\vee b$ and $b\prec a\vee b.$ In the same way, $a\wedge b$ is
the maximal element (with respect to $\prec)$ such that $a\wedge b\prec a$ and
$a\wedge b\prec b$.

\begin{definition}\label{defatom}
 Let $M$ be a monoid. We say that $x\in M$ is an {\em atom} if $x\neq 1$ and if
$x=yz$ implies $y=1$ or $z=1$. $M$ is said to be an {\em atomic monoid} if it is generated
by its atoms and, moreover, for every $a\in M$, there exists an integer $N_a>0$ such that $a$
cannot be written as a product of more than $N_a$ atoms.
\end{definition}

\begin{definition}
 We say that a monoid $M$ is a {\em Gaussian monoid} if it is atomic, (left and right)
cancellative, and if every pair of elements in $M$ admits a (left and right) l.c.m. and
a (left and right) g.c.d.
\end{definition}

\begin{definition}
 A {\em Garside monoid} is a Gaussian monoid which has a {\em Garside element}. A {\em Garside
element} is an element $\Delta\in M$ whose left divisors coincide with their right divisors,
they form a finite set, and they generate $M$.
\end{definition}

\begin{definition}
 The left (and right) divisors of $\Delta$ in a Garside monoid $M$ are called
\textit{simple elements}. We denote by $S$ the (finite) set of simple elements.
\end{definition}

It is known that every Garside monoid admits a group of fractions, and we have:

\begin{definition}
 A group $G$ is called a {\em Garside group} if it is the group of fractions of a Garside monoid.
\end{definition}

The main example of a Garside monoid, as with groups, is the Artin braid
monoid on $n$ strands, $B_{n}^{+}$. It is defined by
Presentation~(\ref{presen}), considered as a presentation for a monoid.
Its group of fractions is the braid group $B_{n}$, and Garside \cite{G} showed that $B_{n}^{+}\subset B_{n}.$
Actually, every Garside monoid embeds into its corresponding Garside group~\cite{DP}.

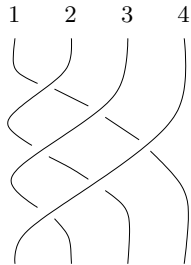
\begin{figure}[ht]
 \centerline{\input{Delta475.pstex_t}}
\caption{The Garside element $\Delta\in B_4^+$.}\label{Delta4}
\end{figure}

Braids in $B_n$ are usually represented as $n$ disjoint strands in ${\mathbb R}^3$, whose endpoints are fixed,
where every horizontal plane between the top and the bottom level intersects each strand in exactly one point,
as in Figure~\ref{Delta4}. Simple elements in $B_n^+$ are easy to recognize: they are those braids
in which any two strands cross at most once. The Garside element of $B_n^+$ is $\Delta =\left( \sigma _{1}\sigma
_{2}\cdots \sigma _{n-1}\right) \left( \sigma _{1}\sigma
_{2}\cdots \sigma _{n-2}\right) \cdots \left( \sigma _{1}\sigma
_{2}\right) \sigma _{1},$ and is represented in
Figure~\ref{Delta4} for $n=4$ (where, as usual, $\sigma_i$
represents a crossing of the strands in positions $i$ and
$i+1$).

There is another important example of Garside monoid, the Birman-Ko-Lee
monoid \cite {BKL}, which has the following presentation:
\begin{equation*}
BKL_{n}^{+}=\left\langle
\begin{array}{ll}
a_{ts} (n\geq t>s\geq 1) & \left|
\begin{array}{l}
 a_{ts}a_{rq}=a_{rq}a_{ts}\text{
if }\left( t-r\right) \left( t-q\right) \left( s-r\right) \left(
s-q\right) >0 \\ a_{ts}a_{sr}=a_{tr}a_{ts}=a_{sr}a_{tr},\text{
where } n\geq t>s>r\geq 1
\end{array}\right.
\end{array}
\right\rangle .  \label{presen2}
\end{equation*}

Its group of fractions is again  the braid group $B_{n}$. The Garside element in
$BKL_{n}^{+}$ is $\delta =a_{n,n-1}a_{n-1,n-2}\cdots a_{2,1}.$ In this monoid
we can perform some computations concerning braid groups faster than using
Artin monoids. Anyway, using the algorithm in~\cite{FGM}, the conjugacy problem has
virtually the same complexity in both monoids.

 From now on, $M$ will denote a Garside monoid, $G$ its group of fractions and $\Delta $ the
corresponding Garside element. Since $M\subset G$, we will refer to the elements in $M$
as the {\em positive} elements of $G$.

%
%

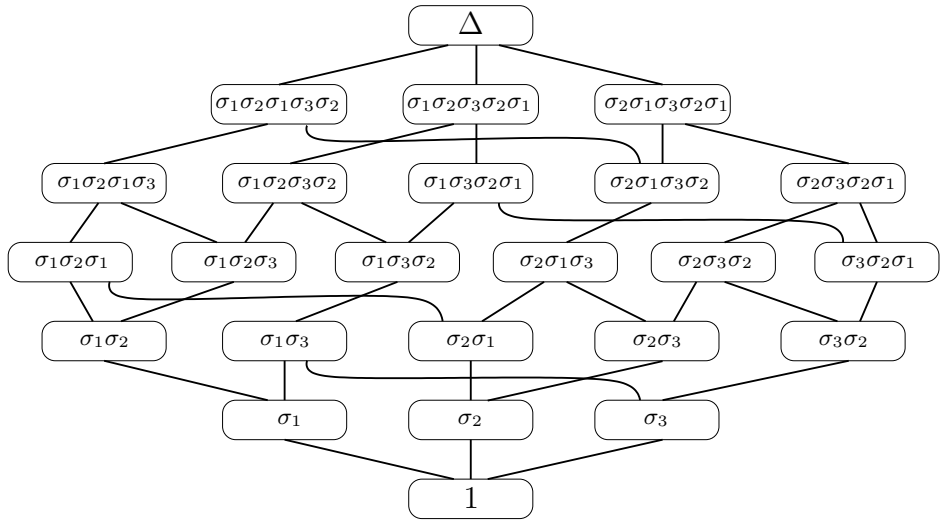
\begin{figure}[ht]
 \centerline{\input{lattice.pstex_t}}
\caption{The lattice of simple elements in
$B_4^+$.}\label{lattice}
\end{figure}

From the existence of l.c.m.'s and g.c.d.'s, it follows that
$\left( M,\prec \right) $ has a lattice structure, and
$S$ becomes a finite sublattice with minimum $1$ and maximum
$\Delta$. In Figure~\ref{lattice} we can see the Hasse diagram of the lattice of simple
elements in $B_{4}^{+}$, where the lines represent left divisibility (from
bottom to top).

We end this section with an important result concerning Garside groups.

\begin{theorem}{\em \cite{DP}}
 For every element $a$ in a Garside group $G$, there exists a unique word
in the atoms of $G$ (and their inverses) representing $a$, called the {\em normal form} of $a$,
and there exists an algorithm that, given a word $w$
in the atoms and their inverses, computes the normal form of the element represented by $w$.
\end{theorem}

\section{Minimal simple elements}

  Simple elements represent a key concept in almost every algorithm concerning Garside groups
(or braid groups): they have been used to compute bi-automatic
normal forms in \cite{T} and \cite{D}, to solve the conjugacy problem in
\cite{Pi2}, \cite{M} and \cite{BKL}, and to compute centralizers
in \cite{Ma}. In some cases, the complexity of these algorithms is
too big due to the size of the set $S$. For instance, in $B_n^+$,
the cardinal of $S$ is $n!$, and this makes the algorithm in
\cite{M} work too slowly. This problem was avoided
in~\cite{FGM}, by considering {\em minimal simple elements}. We
will also use minimal simple elements in this paper, so this
section is devoted to them.

 Given an element $a$ in a Garside group $G$, there exists a subset $C^{sum}(a)$
of the conjugacy class of $a$, called {\em Summit Class of $a$}, satisfying some suitable properties.
In~\cite{M}, when talking about braids, this subset is called {\em Super Summit Set}, but when we talk
about Garside groups we prefer to use the terminology in~\cite{Pi2}. Roughly speaking, $C^{sum}(a)$ is
the set of conjugates of $a$ having the `simplest' normal form, in a certain sense. Hence,
$C^{sum}(a)$ is an invariant of the conjugacy class of $a$ (it does not depend on $a$, but on its
conjugacy class).

 There exists a procedure, called `cycling and decycling', to obtain an element $a'\in C^{sum}(a)$ and
an element $x$ such that $x^{-1}ax=a'$ (see~\cite{M}). The centralizers of $a$ and $a'$ are then
related as follows: $Z(a)=x Z(a') x^{-1}$. Hence, if we know a generating set for $Z(a')$, we obtain
immediately a generating set for $Z(a)$, with the same number of elements: it suffices to conjugate every
generator by $x$. Therefore, we will just study the elements in the Summit Class of $a$.

 Consider an element $v\in C^{sum}(a)$. If we conjugate $v$ by a nontrivial simple element, we obtain
an element in $G$, that may or may not be in $C^{sum}(a)$. We will consider just the elements in
$S\backslash\{1\}$ that conjugate $v$ to an element in $C^{sum}(a)$. Among these simple elements, we take
those which are minimal with respect to $\prec$, and we call this set $S_v^{sum}$.
In other words, we define $S_v^{sum}$ as the set of minimal elements (with respect to $\prec$) in
$\left\{ s\in S\backslash\{1\}:\; s^{-1}vs\in C^{sum}(a)\right\}$.

 There are two important results concerning these minimal simple elements:

\begin{proposition}{\em \cite{FGM}}
 Let $M$ be a Garside monoid with $t$ atoms, $G$ its corresponding Garside group, and $a\in G$.
For every $v\in C^{sum}(a)$, the cardinal of $S_v^{sum}$ is no bigger than $t$.
\end{proposition}

\begin{proposition}{\em \cite{FGM}}\label{prochain}
 Let $u,v$ be two conjugate elements in $C^{sum}(a)$. Then there exists a sequence
$u=u_{1}, u_2,...,u_{k}=v$ of elements in $C^{sum}(a)$ such that, for
$i=1,...,k-1$, there exists $s_{i}\in S_{u_{i}}^{sum}$ verifying
$u_{i}s_{i}=s_{i}u_{i+1}.$
\end{proposition}

 We will represent the above property as follows:
$$
  u=u_1 \stackrel{s_1}{\longrightarrow} u_2 \stackrel{s_2}{\longrightarrow} u_3 \rightarrow \cdots \rightarrow
  u_{k-1} \stackrel{s_{k-1}}{\longrightarrow} u_k=v,
$$
where $s_i\in S_{u_i}^{sum}$ for every $i$, and the arrow means conjugation by the corresponding $s_i$.
We call such a sequence a {\em minimal chain} from $u$ to $v$.

\begin{example}
 Consider the braid monoid $B_4^+$. As we saw in the previous section, the set of simple elements in
$B_4^+$ has $24$ elements (see Figure~\ref{lattice}). Consider $\sigma_1\in C^{sum}(\sigma_1)\subset B_4$.
Then $S_{\sigma_1}^{sum}=\{\sigma_1, \: \sigma_2\sigma_1, \: \sigma_3\}$. The conjugates of $\sigma_1$ by these
three elements are, respectively, $\sigma_1$, $\sigma_2$ and $\sigma_1$. All of them lie in
$C^{sum}(\sigma_1)$. The conjugating elements are clearly minimal: $\sigma_1$ and $\sigma_3$ do not have
nontrivial divisors, and the only nontrivial divisor of $\sigma_2 \sigma_1$ is $\sigma_2$, which does not
conjugate $\sigma_1$ to a positive element (hence to an element in $C^{sum}(\sigma_1)$).
\end{example}

\begin{remark}
 It is shown in~{\em \cite{FGM}} that for every $v\in C^{sum}(a)$ and every atom $x$, there exists at
most one element $s\in S_v^{sum}$ which is a multiple of $x$. This is why the cardinal of $S_v^{sum}$
is bounded by the number of atoms. In $B_n^+$, the atoms are $\sigma_1,\ldots,\sigma_{n-1}$, and in the
above example we can clearly see which element in $S_{\sigma_1}^{sum}$ corresponds to each atom.
\end{remark}

 In general, for a given $v\in C^{sum}(a)$, there are strictly fewer minimal simple elements than atoms, as we
can see in the following example:

\begin{example}
  Let $v=\sigma_1 \sigma_2\in C^{sum}(\sigma_1 \sigma_2)\subset B_4^+$.
Then $S_v^{sum}=\{\sigma_1, \sigma_3\sigma_2\sigma_1\}$. Indeed, conjugating we obtain
$\sigma_1^{-1} (\sigma_1\sigma_2) \sigma_1= \sigma_2\sigma_1$, and
$(\sigma_3\sigma_2\sigma_1)^{-1}\sigma_1\sigma_2(\sigma_3\sigma_2\sigma_1)=\sigma_2\sigma_3$.
But the minimal multiple of $\sigma_2$ which conjugates $v$ to a positive element is
$\sigma_2\sigma_1\sigma_2$, which is also a multiple of $\sigma_1$, so it is not in $S_v^{sum}$
(since it is not minimal).
\end{example}

 In \cite{FGM} the authors give an algorithm to compute $S_v^{sum}$, given $v\in C^{sum}(a)$, and use it to
compute the whole Summit Class of any element. Sometimes, a problem can be solved using either simple
elements, or minimal simple elements. The latter possibility is usually much faster. For instance, in the
braid monoid $B_n^+$, computing the set $S_v^{sum}$ takes time $O(l^2n^4)$, where $l$ is the word-length of
$v$. After performing this fast computation, we can work with a set of less than $n-1$ elements
($S_v^{sum}$), instead of a set with $n!$ elements ($S$).

 In order to compute centralizers in Garside groups, Makanin~\cite{Ma} used simple elements,
but we are going to see in the next section how the use of minimal simple elements, and a new approach
to the problem, can make the computations much faster.

\section{Minimal summit graph}

We shall explain in this section a new approach to our problem,
which involves the fundamental group of a certain graph. Consider
an element $a$ in a Garside group $G$. We want to find a
generating set for the centralizer of $a$. As we said in
Section~\ref{secsimple}, we will study the elements in its
Summit Class $C^{sum}(a)$.

 Let us construct a directed graph $\Gamma$, that we call {\em minimal summit graph} of $a$. The vertices
of $\Gamma$ are the elements in $C^{sum}(a)$. The arrows of $\Gamma$
are labelled by simple elements, in the following way: For every two vertices
$v$ and $w$, an arrow labelled by $s$ goes from $v$ to $w$ if and only if
$s\in S_v^{sum}$ and $s^{-1}vs=w$. In other words, $s$ is a minimal simple element that conjugates $v$ to an
element in $C^{sum}(a)$, and $w$ is the result of that conjugation.
Therefore, every path in $\Gamma$ going from a vertex $u$ to another vertex
$v$, and moving always in the sense of the arrows, is a minimal chain from $u$ to $v$ (see
Proposition~\ref{prochain}).

\begin{figure}[ht]
\centerline{\input{MCG1.pstex_t}}
\caption{Minimal summit graph of $\sigma_1\in B_4^+$.}
\label{MCG1}
\end{figure}
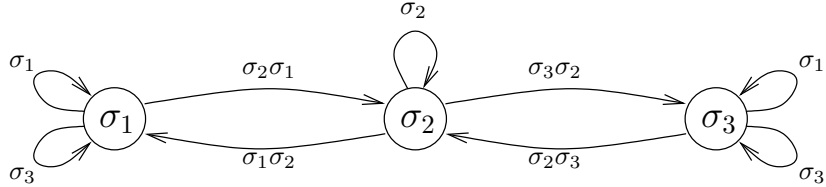

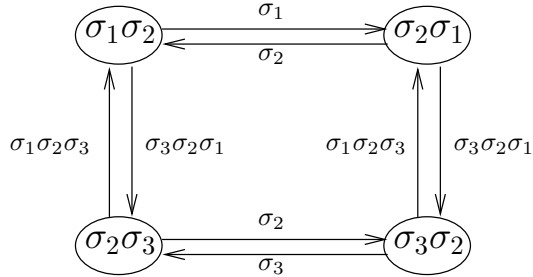
\begin{figure}[ht]
\centerline{\input{MCG12.pstex_t}}
\caption{Minimal summit graph of $\sigma_1\sigma_2\in B_4^+$.}
\label{MCG12}
\end{figure}

The minimal summit graph of $\sigma_1\in B_4^+$ is represented in Figure~\ref{MCG1}, and that of
$\sigma_1\sigma_2$ in Figure~\ref{MCG12}.

The main idea in our algorithm is the following: Given $a'\in C^{sum}(a)$, every element in
$Z(a')$ can be seen as a loop in $\Gamma$, based at $a'$. So every
generating set for the fundamental group of $\Gamma$ corresponds
to a generating set for $Z(a')$ (recall that if we know a generating set for $Z(a')$, we
also know a generating set for $Z(a)$). We devote the rest of this section to proving this.
We shall need the following results:

\begin{lemma}\label{lempositgen}
 For every $a\in G$, the centralizer of $a$ can be generated by elements in $M$.
\end{lemma}

\begin{proof}
 Let $c\in Z(a)$. We will try to write $c$ as a product of positive elements in $Z(a)$ (and their inverses).
We know by~\cite{DP} that there is an integer $k$ such that $\Delta^k$ is in
the center of $G$ (thus in $Z(a)$), and another integer $r$, big
enough, such that $\Delta^{kr}c\in M$. Hence,
$c=(\Delta^{kr})^{-1}(\Delta^{kr}c)$, where $\Delta^{kr}$ and
$\Delta^{kr}c$ belong to $M\cap Z(a)$. This implies the result.
\end{proof}

\begin{theorem}{\em \cite{Pi2}}\label{teomax}
  Let $u,v\in C^{sum}(a)$ and $x\in M$ such that $x^{-1}ux=v$. Let $s\in S$ be the maximal simple prefix
of $x$, that is, $s$ is maximal (with respect to $\prec$) among the simple elements dividing $x$. Then
$s^{-1}us\in C^{sum}(a)$.
\end{theorem}

\begin{corollary}\label{corchain}
 Let $u,v \in C^{sum}(a)$, and $x\in M$ as above. Then there exists a decomposition $x=s_1s_2\cdots s_{k-1}$,
and $k$ elements $u=u_1,u_2,\ldots,u_k=v\in C^{sum}(a)$, such that
$$
  u=u_1 \stackrel{s_1}{\longrightarrow} u_2 \stackrel{s_2}{\longrightarrow} u_3 \rightarrow \cdots \rightarrow
  u_{k-1} \stackrel{s_{k-1}}{\longrightarrow} u_k=v,
$$
is a minimal chain from $u$ to $v$.
\end{corollary}

\begin{proof}
  First, let us decompose $x=t_1t_2\cdots t_{p-1}$, where for every $i$, $t_i$ is the maximal simple prefix
of $t_it_{i+1}\cdots t_{p-1}$ (this is the {\em left greedy normal form} of $x$, in the sense of~\cite{T}).
By Theorem~\ref{teomax}, we obtain a chain
$$
  u=w_1 \stackrel{t_1}{\longrightarrow} w_2 \stackrel{t_2}{\longrightarrow} w_3 \rightarrow \cdots \rightarrow
  w_{p-1} \stackrel{t_{p-1}}{\longrightarrow} w_p=v,
$$
where $w_i\in C^{sum}(a)$ for $i=1,\ldots,p$. But this chain is not necessarily minimal. Now, for every $t_i$, we
proceed as follows: if it is minimal (among the simple elements that conjugate $w_i$ to an element in
$C^{sum}(a)$), we do not touch it. Otherwise, there exists an element $r\in S_{w_i}^{sum}$ dividing $t_i$.
So we can decompose the arrow $w_i \stackrel{t_i}{\longrightarrow} w_{i+1}$ as
$w_i \stackrel{r}{\longrightarrow} w' \stackrel{r'}{\longrightarrow} w_{i+1}$, where $t_i=r\:r'$ and
$w'\in C^{sum}(a)$. If $r'$ is not minimal, we
decompose it in the same way. If we continue this process we obtain, at each step, a decomposition
$t_i=r_1\cdots r_m$, where every $r_j$ is a simple element. Hence we have a chain
$r_1\prec r_1r_2 \prec r_1r_2r_3\prec \cdots \prec (r_1\cdots r_m)$ of simple elements. But the length
of such a chain is bounded above, since there is a finite number of simple elements. Therefore, we cannot
decompose $t_i$ indefinitely, and this process must stop.

 At the end, we will have decomposed every $t_i$ as a product of minimal simple elements, so the
result follows.
\end{proof}

We can finally prove the main result of this section.  Consider the natural group homomorphism
$p:\:\pi_1(\Gamma,a')\longrightarrow G$, which sends every loop in $\Gamma$ based at $a'$ to the element
in $G$ obtained by reading the labels in the path, with the corresponding signs. One has the following:


\begin{theorem}
 The homomorphism $p$ maps $\pi_1(\Gamma,a')$ onto $Z(a')$.
\end{theorem}

\begin{proof}
 Since every loop $\gamma\in \pi_i(\Gamma,a')$ starts and ends at $a'$, then
$p(\gamma)$ conjugates $a'$ to itself, so $p(\gamma)\in Z(a')$.
Hence, we get $p:\:\pi_1(\Gamma,a')\longrightarrow Z(a')$.

  By Lemma~\ref{lempositgen}, we know that $Z(a')$ is generated by positive elements. Every positive
element $y\in Z(a')$ verifies $y^{-1}a'y=a'$, so by
Corollary~\ref{corchain}, $y$ can be decomposed into minimal
simple elements, yielding a minimal chain from $a'$ to itself. This
minimal chain is actually an element in $p^{-1}(y)$. Hence, there
exist preimages by $p$ for all positive elements in $Z(a')$. Since
the positive elements generate $Z(a')$, we get that $p$ is a
surjection, and we are done.
%
\end{proof}

By the above result, in order to compute a generating set for
$Z(a')$ we just need to compute a generating set for
$\pi_1(\Gamma,a')$. It is well known how to do this (see, for
instance, \cite{LS}): Choose a maximal tree $T$ in $\Gamma$. For
every vertex $v$ in $\Gamma$, call $\gamma_v$ the only simple path
in $T$ going from $a'$ to $v$. Let $A$ be the set of arrows in
$\Gamma\backslash T$ and, for every $\alpha\in A$, denote $s(\alpha)$ and
$t(\alpha)$ the {\em starting vertex} and the {\em target} of $\alpha$,
respectively. Then there is a generating set $F$ for
$\pi_1(\Gamma,a')$, which is in one-to-one correspondance with $A$.
It is the following: $F=\{\gamma_{s(\alpha)}\:\alpha\:\gamma_{t(\alpha)}^{-1};\;
\alpha\in A\}$. So $p(F)$ is the generating set for $Z(a')$ that our
algorithm will compute.

\begin{remark}
 In a previous version of this paper, we considered the whole conjugacy class of $a$ (in $M$) instead
of its Summit Class, hence we computed the {\em minimal conjugacy graph} instead of
the minimal summit graph. Although there are no known bounds for the sizes of these sets, the Summit Class
is in general much smaller than the whole conjugacy class, so this new approach is more efficient.
We thank A. Kalka for his observation on this matter.
\end{remark}

\section{The algorithm}

  We shall now explain our algorithm in detail. Let $a$ be an element of a Garside group $G$, and
let $a'\in C^{sum}(a)$. Let $\Gamma$ be the minimal summit graph of $a'$. We will start by
computing $\Gamma$ and, for every vertex $v\in \Gamma$, a path $\gamma_v$ going from $a'$ to $v$ in
a maximal tree $T$ in $\Gamma$.

 In the following routine, $v$ denotes the current vertex of $\Gamma$ under study, $U$ is the set
of known vertices of $\Gamma$ (i.e. the known elements in $C^{sum}(a)$), and
$V$ is the set of vertices which have already been studied by the routine.

\vspace{.3cm}
\noindent
{\bf Routine 1: Computation of $\mathbf{\Gamma}$ and $\mathbf{T}$.}

\vspace{.3cm}
Input: $a'\in C^{sum}(a)$.

\begin{enumerate}

 \item Set $v=a'$, $U=\{a'\}$, $V=\phi$, $\Gamma=\phi$, $T=\phi$ and $\gamma_a'=1$.

 \item Compute $S_v^{sum}$.

 \item For every $s\in S_v^{sum}$ do the following:
 \begin{enumerate}

    \item Set $w=s^{-1}vs\in C^{sum}(a)$, written in normal form. Set $\Gamma=\Gamma\cup \{(v,s,w)\}$.

    \item If $w\notin U$, set $U=U \cup \{w\}$, $T=T\cup \{(v,s,w)\}$ and $\gamma_w=\gamma_vs$.

 \end{enumerate}

 \item Set $V=V\cup \{v\}$. If $U\neq V$, take an element $x\in U\backslash V$. Set $v=x$ and go to
       Step 2.

 \item Stop.

\end{enumerate}

 From the results in the previous sections we can see that, at the end of this routine, we will obtain
the following data:

\begin{itemize}

 \item A set $U=V=C^{sum}(a)$, which is the set of vertices of $\Gamma$.

 \item A set $\Gamma$ which corresponds to the graph $\Gamma$: it
contains an element $(v,s,w)$ for every arrow of the graph $\Gamma$ labelled
by $s$, and going from $v$ to $w$.

 \item A set $T$ which corresponds to a subgraph of $\Gamma$.

 \item For every $v\in V$, a path $\gamma_v$ in the subgraph $T$, going from $a'$ to $v$.

\end{itemize}

\begin{proposition}
 $T$ is a maximal tree in $\Gamma$.
\end{proposition}

\begin{proof}
 The graph $T$ is computed by Routine~1 as follows: Let $s$ be an arrow such that $t(s)=w\neq a'$.
Then $s$ is added to $T$ (in Step 3(b)) if and only if it is the first arrow considered by Routine~1
whose target is $w$. Hence, for every $w\in V$, $w\neq a'$, there is exactly one arrow in $T$ ending
at $w$. And there is no arrow in $T$ ending at $a'$.

Therefore, if we start at a vertex $v$, and we try to construct a path in $T$, as long as possible, moving
always in the sense opposite to the arrows, we have a unique choice. This path would always end at $a'$,
and it is actually the inverse of the path $\gamma_v$ computed by Routine~1: Just notice that the
path $\gamma_v$ goes from $a'$ to $v$ always in the sense of the arrows.

 Let us then show that $T$ is a tree. Suppose that there exists a nontrivial simple loop $\alpha$ in $T$.
Since there is no pair of arrows of $T$ with the same target, we can assume that $\alpha$ moves always in
the sense of the arrows. Since $a'$ is not the
target of any arrow in $T$, then $a'$ does not belong to the set of vertices in $\alpha$. But, if we start
at a vertex $v$ in $\alpha$, and we try to follow $\gamma_v^{-1}$ as above, we would go along $\alpha^{-1}$
an infinite number of times, never reaching $a'$. This contradiction shows that there are no loops in
$T$, so it is a tree.

Finally, $T$ is maximal since it is connected (every vertex is connected to $a'$), and it contains all the
vertices in $\Gamma$.
\end{proof}

Therefore, we can use the data given by Routine~1 to compute a
generating system for $Z(a)$, by the procedure explained in the
previous section:

\vspace{.3cm} \noindent {\bf Routine 2: Computation of a generating set for
$\mathbf{Z(a)}$.}

\vspace{.3cm}
Input: $a\in G$.

\begin{enumerate}

 \item Using `cyclings and decyclings', compute $a'\in C^{sum}(a)$, and $x\in G$ such that
     $x^{-1} a x = a'$.

 \item Apply Routine~1 to $a'$, obtaining $\Gamma$, $T$ and the paths $\gamma_v$.

 \item Set $N=\phi$.

 \item For every $(v,s,w) \in \Gamma\backslash T$ do the following:

 \begin{enumerate}

   \item  Compute the normal form $\alpha$ of $x (\gamma_v s \gamma_w^{-1}) x^{-1}$
         (given as an element of $G$).

   \item  If $\alpha\notin N$, set $N=N\cup \{\alpha\}$.

 \end{enumerate}

 \item Return $N$. Stop.

\end{enumerate}

\section{Complexity}

 In order to study the complexity of our algorithm, we should know some data concerning the Garside
monoid $M$, and the element $a\in G$ under study:

\begin{itemize}

 \item $t$: The number of atoms in $M$.

 \item $m$: The maximal length of a simple element in $M$.

 \item $k$: The number of elements in $C^{sum}(a)$ (i.e. the number of vertices in $\Gamma$).

 \item $l$: The maximal word length of an element in $C^{sum}(a)$.

 \item $D$: The complexity of computing $a'\in C^{sum}(a)$ and $x$, as in Routine~2.

 \item $C$: The complexity of an algorithm to compute $S_v^{sum}$ for an element $v$ of word length $l$.

 \item $N_i$: The complexity of computing the normal form of a word of length $i$.

\end{itemize}

 If we know all the previous data, we can compute the complexity of our algorithm by the following result:

\begin{proposition}
 Given an element $a$ in a Garside group $G$, we can compute a generating set for $Z(a)$ in time
$O(D+kC+ktN_{2km})$.
\end{proposition}

\begin{proof}
We start by computing $a'$ and $x$, taking time $O(D)$. Then we run Routine~1, which does the following:
For every vertex $v$ in $T$, it computes $S_v^{sum}$ and then, for every $s\in S_v^{sum}$, it computes
the normal form of $s^{-1}vs$. The other steps in Routine~1 are negligible. Moreover, the algorithm
used in~\cite{FGM} to compute $S_v^{sum}$ also gives the normal forms of $s^{-1}vs$, for $s\in S_v^{sum}$.
Hence, Routine~1 has complexity $O(kC)$.

 Routine~2 continues by computing the normal form of $x(\gamma_v s \gamma_w^{-1})x^{-1}$,
for every arrow $(v,s,w)$ in
$\Gamma \backslash T$. We know that the number of arrows in $\Gamma$ is bounded by $kt$, since there
are at most $t$ arrows for each vertex, and there are $k$ vertices. On the other hand, there is exactly
one arrow in $T$ whose target is $v$, for every vertex in $\Gamma$ different from $a'$. Hence, $T$ has
$k-1$ arrows, so $\Gamma\backslash T$ has at most $kt-k-1$ arrows. Now $x(\gamma_v s \gamma_w^{-1})x^{-1}$
is a product of at most $2(k+|x|)-1$ simple elements. Hence, written as a word in the atoms and their
inverses, its length is bounded by $2(k+|x|)m$. Since the length of $x$ is negligible compared to $k$,
$N_{2(k+|x|)m}$ is equivalent to $N_{2km}$. Therefore, the complexity of this loop is $O(ktN_{2km})$,
and the result follows.
\end{proof}

 For some particular Garside monoids and groups, one would like to know the complexity in more detail,
just depending on the word length of $a$, and on some integer related to the monoid. This can be done
more easily for Garside monoids in which every relation is homogeneous, for in this case, all the elements
in the conjugacy class of $a$ have the same word length. This is the case for $B_n^+$, $BKL_n^+$ and
Artin monoids. The authors have studied in~\cite{FGM} the complexity $C$ of computing $S_v^{sum}$ for a
braid $v$ of length $l$ (either in $B_n^+$ or in $BKL_n^+$). In~\cite{T} and in \cite{BKL} we can find the
complexity $N_i$, for elements in $B_n^+$ and in $BKL_n^+$ respectively, and in~\cite{BKL2}
the complexity $D$ is given. Hence, we obtain the following
results:

\begin{corollary}
  Given $a\in B_n$ of word length $l$ in the Artin generators, the complexity of computing $Z(a)$
(using the Garside structure given by $B_n^+$) is $O(k^3l^2n^6\log n)$.
\end{corollary}

\begin{proof}
  In $B_n^+$, one has $t=n-1$, $m=\frac{n(n-1)}{2}$, $C=O(l^2n^4)$ (see~\cite{FGM}),
$N_i=O(i^2n\log n)$ (see~\cite{T}) and $D=O(l^2n^3)$ (see~\cite{BKL2}). Hence, the complexity of our
algorithm to compute $Z(a)$ becomes $O\left(l^2n^3+kl^2n^4+k(n-1)(kn(n-1))^2n\log n\right) =
O(kl^2n^4+k^3n^6\log n) = O(k^3l^2n^6\log n)$, so the result is true.
\end{proof}

\begin{corollary}
  Given $a\in B_n$ of word length $l$ in the Birman-Ko-Lee generators, the complexity of computing
$Z(a)$ (using the Garside structure given by $BKL_n^+$) is $O(k^3l^2n^5)$.
\end{corollary}

\begin{proof}
  This time, in $BKL_n^+$, one has $t=\frac{n(n-1)}{2}$, $m=n-1$, $C=O(l^2n^5)$ (see~\cite{FGM}),
$N_i=O(i^2n)$ (see~\cite{BKL}) and $D=O(l^2n^2)$ (see~\cite{BKL2}). Therefore, our algorithm for
computing $Z(a)$ has complexity $O\left(l^2n^2+ kl^2n^5 +k\frac{n(n-1)}{2}(2k(n-1))^2n\right) =
O(kl^2n^5 +k^3n^5) = O(k^3l^2n^5)$.
\end{proof}

It would remain to know, in both cases, a bound for $k$ in terms of $l$ and $n$. This is still not
known, but Thurston, in~\cite{T}, conjectures that $k$ is bounded by a polynomial in $l$ (although it
seems to be exponential in $n$).

\section{Effective computations}

 In this section we show the results we have obtained by implementing our algorithm. We have computed
generating sets for the centralizers of many positive elements in
the braid monoids $B_n^+$, for $n=3,\ldots,8$. We have been
exhaustive, computing centralizers of all braids of a given
length, in order to conjecture an upper bound for the number of
generators.

 We proceeded as follows: first, by using the algorithm in~\cite{FGM}, we computed the conjugacy
classes in $B_n^+$ of all braids of the considered length. Notice
that two elements in the same conjugacy class have conjugated
centralizers: if $c^{-1}ac=b$ and $x\in Z(a)$, then $c^{-1}xc\in
Z(b)$; hence, the number of generators in the centralizer of $a$
and $b$ are the same. Therefore, we just had to compute the
centralizer of one representative for each conjugacy class in
$B_n^+$.

The results of these computations were surprising, since the
number of generators were quite small. In the following table we
can see the braids that we tested, and the maximal size of a
generating set for the centralizer, in each case:

\vspace{.3cm}
\begin{tabular}{|c|c|c|c|}
\hline
n & Length of braids & Number of Conj. Classes & Max. number of generators \\
\hline
\hline
 3 & $4\leq l \leq 20$ & 1634 &  4  \\
\hline
 4 & $4 \leq l \leq 15$ & 4225 &  16 \\
\hline
 5 & $4 \leq l \leq 12 $ & 2314 &  17 \\
\hline
 6 &  $4 \leq l \leq 10$ & 1152 &  12  \\
\hline
 7 &  $4 \leq l \leq 10$ & 1753 &  17  \\
\hline
 8 &  $4 \leq l \leq 8$ & 521 & 22  \\
\hline
\end{tabular}

\vspace{.5cm}
 Actually, we found out that the generators obtained by the algorithm were not always independent, so we were able to eliminate
some of them. For instance, if we compute $Z(a)$ for
$a=\sigma_1\in B_4$, the algorithm will give the following
generating set:
$$
\{\sigma_1, \;\sigma_2\sigma_1\sigma_1\sigma_2, \;\sigma_3,
\;\sigma_2\sigma_1(\sigma_3\sigma_2\sigma_2\sigma_3)\sigma_1^{-1}\sigma_2^{-1}\}.
$$

But the fourth element can also be written as $(\sigma_3)^{-1}(\sigma_2\sigma_1\sigma_1\sigma_2)(\sigma_3)$,
so it can be eliminated from the generating set, yielding:
$$
Z(\sigma_1)=\left<\sigma_1, \;\sigma_2\sigma_1\sigma_1\sigma_2,
\;\sigma_3\right>\subset B_4.
$$

In the case of $B_3$, we were able to obtain the following: for
every positive braid $a\in B_3^+$ of length $l\leq 20$, there is a
generating set for $Z(a)$ with at most two elements.

 We cannot show here all the results but we can see, as an example, the following table: it contains a
representative for each conjugacy class of elements in $B_3^+$ of
length 11, and a generating set for their centralizers.

{\small
\vspace{.3cm}
\begin{tabular}{|c||c|c|c|}
\hline
\multicolumn{3}{|c|}{Centralizers of braids in $B_3^+$ of length 11} \\
\hline
\hline
$a$ & \multicolumn{2}{|c|}{Generators for $Z(a)$} \\
\hline
\hline
$\sigma_1^{11}$ & $\sigma_1$ & $\sigma_2\sigma_1^{2}\sigma_2$   \\
\hline
$\sigma_1^{10}\sigma_2$ & $\sigma_1\sigma_2\sigma_1^{2}\sigma_2\sigma_1$ & $\sigma_1^{2}\sigma_2^{2}\sigma_1^{2}\sigma_2\sigma_1^{-6}$   \\
\hline
$\sigma_1^{9}\sigma_2^{2}$ & $\sigma_1\sigma_2\sigma_1^{2}\sigma_2\sigma_1$ & $\sigma_1^{6}\sigma_2^{-1}\sigma_1^{-2}\sigma_2^{-3}\sigma_1^{-1}$   \\
\hline
$\sigma_1^{8}\sigma_2^{3}$ & $\sigma_1\sigma_2\sigma_1^{2}\sigma_2\sigma_1$ & $\sigma_1^{6}\sigma_2^{-1}\sigma_1^{-3}\sigma_2^{-2}\sigma_1^{-1}$   \\
\hline
$\sigma_1^{2}\sigma_2^{6}\sigma_1^{2}\sigma_2$ & $\sigma_1^{2}\sigma_2\sigma_1^{-2}$ & $\sigma_1^{3}\sigma_2^{2}\sigma_1^{-1}$   \\
\hline
$\sigma_1^{7}\sigma_2^{4}$ & $\sigma_1\sigma_2\sigma_1^{2}\sigma_2\sigma_1$ & $\sigma_1^{6}\sigma_2^{-2}\sigma_1^{-2}\sigma_2^{-2}\sigma_1^{-1}$   \\
\hline
$\sigma_1^{6}\sigma_2^{2}\sigma_1^{2}\sigma_2$ & $\sigma_1^{4}\sigma_2\sigma_1^{-4}$ & $\sigma_1\sigma_2\sigma_1^{2}\sigma_2\sigma_1$   \\
\hline
$\sigma_1^{6}\sigma_2^{5}$ & $\sigma_1\sigma_2\sigma_1^{2}\sigma_2\sigma_1$ & $\sigma_1^{6}\sigma_2\sigma_1^{-1}\sigma_2^{-2}\sigma_1^{-2}\sigma_2^{-2}\sigma_1^{-1}$   \\
\hline
$\sigma_1^{5}\sigma_2^{3}\sigma_1^{2}\sigma_2$ & $\sigma_1^{3}\sigma_2^{2}\sigma_1^{-4}$ & $\sigma_1\sigma_2\sigma_1^{2}\sigma_2\sigma_1$   \\
\hline
$\sigma_1^{5}\sigma_2^{2}\sigma_1^{3}\sigma_2$ & $\sigma_1^{2}\sigma_2^{3}\sigma_1^{-4}$ & $\sigma_1\sigma_2\sigma_1^{2}\sigma_2\sigma_1$   \\
\hline
$\sigma_1^{5}\sigma_2^{2}\sigma_1^{2}\sigma_2^{2}$ & $\sigma_1\sigma_2\sigma_1^{2}\sigma_2\sigma_1$ & $\sigma_1\sigma_2^{4}\sigma_1^{-4}$   \\
\hline
$\sigma_1^{4}\sigma_2^{2}\sigma_1^{4}\sigma_2$ & $\sigma_1\sigma_2\sigma_1^{2}\sigma_2\sigma_1$ & $\sigma_1^{2}\sigma_2^{2}\sigma_1^{2}\sigma_2^{-1}\sigma_1^{-4}$   \\
\hline
$\sigma_1^{4}\sigma_2^{2}\sigma_1^{3}\sigma_2^{2}$ & $\sigma_1\sigma_2\sigma_1^{2}\sigma_2\sigma_1$ & $\sigma_1\sigma_2^{3}\sigma_1^{2}\sigma_2^{-1}\sigma_1^{-4}$   \\
\hline
$\sigma_1^{4}\sigma_2^{2}\sigma_1^{2}\sigma_2^{3}$ & $\sigma_1\sigma_2\sigma_1^{2}\sigma_2\sigma_1$ & $\sigma_1\sigma_2^{2}\sigma_1^{3}\sigma_2^{-1}\sigma_1^{-4}$   \\
\hline
$\sigma_1^{4}\sigma_2^{3}\sigma_1^{2}\sigma_2^{2}$ & $\sigma_1\sigma_2\sigma_1^{2}\sigma_2\sigma_1$ & $\sigma_1\sigma_2^{4}\sigma_1\sigma_2^{-1}\sigma_1^{-4}$   \\
\hline
$\sigma_1^{3}\sigma_2^{2}\sigma_1^{3}\sigma_2^{3}$ & $\sigma_1\sigma_2\sigma_1^{2}\sigma_2\sigma_1$ & $\sigma_1^{3}\sigma_2^{2}\sigma_1\sigma_2^{-1}\sigma_1^{-3}\sigma_2^{-2}\sigma_1^{-1}$   \\
\hline
\end{tabular}
}

\vspace{.5cm}
 When $n$ becomes bigger, it is more difficult to eliminate generators by hand. Nevertheless, we can show as an example the
following table, where we can see a representative for every conjugacy class of elements of length $6$ in $B_4^+$. We were
able to reduce the number of generators to be less than or equal to $3$ in every case:

\vspace{.5cm}
{\small
\begin{tabular}{|c||c|c|c|}
\hline
\multicolumn{4}{|c|}{Centralizers of braids in $B_4^+$ of length 6} \\
\hline
\hline
$a$ & \multicolumn{3}{|c|}{Generators for $Z(a)$} \\
\hline
\hline
$\sigma_1^{6}$ & $\sigma_1$ & $\sigma_3$ & $\sigma_2\sigma_1^{2}\sigma_2$  \\
\hline
$\sigma_1^{5}\sigma_2$ & $\sigma_3\sigma_2\sigma_1^{2}\sigma_2\sigma_3$ & $\sigma_1^{2}\sigma_2\sigma_1^{-3}$ & $\sigma_1^5 \sigma_2$ \\
\hline
$\sigma_1^{5}\sigma_3$ & $\sigma_1$ & $\sigma_3$ & $\sigma_2\sigma_1\sigma_3\sigma_2^{2}\sigma_1\sigma_3\sigma_2$   \\
\hline
$\sigma_1^{4}\sigma_2^{2}$ & $\sigma_3\sigma_2\sigma_1^{2}\sigma_2\sigma_3$ & $\sigma_1\sigma_2\sigma_1^{2}\sigma_2\sigma_1$ & $\sigma_1^{4}\sigma_2^{2}$   \\
\hline
$\sigma_1^{4}\sigma_2\sigma_3$ & $\sigma_1^{2}\sigma_2\sigma_1\sigma_3\sigma_2^{-2}\sigma_1^{-3}$ & $\sigma_1^{4}\sigma_2\sigma_3$ & $\sigma_1\sigma_2\sigma_1^{2}\sigma_2\sigma_1\sigma_3\sigma_2\sigma_1^{-2}$   \\
\hline
$\sigma_1^{3}\sigma_3\sigma_1\sigma_3$ & $\sigma_1$ & $\sigma_3$ & $\sigma_2\sigma_1\sigma_3\sigma_2^{2}\sigma_1\sigma_3\sigma_2$   \\
\hline
$\sigma_1^{3}\sigma_2^{3}$ & $\sigma_1\sigma_2\sigma_1^{-2}$ & $\sigma_3\sigma_2\sigma_1^{2}\sigma_2\sigma_3$ & $\sigma_1\sigma_2\sigma_1^{2}\sigma_2\sigma_1$   \\
\hline
$\sigma_1^{3}\sigma_2^{2}\sigma_3$ & $\sigma_1\sigma_2\sigma_1\sigma_3\sigma_1\sigma_2\sigma_3\sigma_2\sigma_1^{-2}$ & $\sigma_1^{3}\sigma_2^{2}\sigma_3$ & \\
\hline
$\sigma_1^{3}\sigma_2\sigma_3\sigma_2$ & $\sigma_1^{2}\sigma_3\sigma_2\sigma_3^{-1}\sigma_2^{-2}\sigma_1^{-1}$ & $\sigma_1\sigma_2\sigma_1\sigma_2\sigma_1\sigma_3\sigma_2\sigma_1^{-1}$ & $\sigma_1^3 \sigma_2\sigma_3\sigma_2$\\
\hline
$\sigma_1\sigma_3\sigma_1\sigma_3\sigma_1\sigma_3$ & $\sigma_1$ & $\sigma_2\sigma_1\sigma_3\sigma_2$ & $\sigma_3$   \\
\hline
$\sigma_1\sigma_2\sigma_1^{2}\sigma_2\sigma_1$ & $\sigma_1$ & $\sigma_2$ & $\sigma_3\sigma_2\sigma_1^{2}\sigma_2\sigma_3$   \\
\hline
$\sigma_1^{2}\sigma_2\sigma_1\sigma_3\sigma_2$ & $\sigma_1\sigma_2\sigma_1^{-1}$ & $\sigma_1\sigma_3$ & $\sigma_2\sigma_3\sigma_2^{-1}$  \\
\hline
$\sigma_1^{2}\sigma_2^{3}\sigma_3$ & $\sigma_1\sigma_2\sigma_1\sigma_3\sigma_1\sigma_2\sigma_3\sigma_2\sigma_1\sigma_2^{-1}\sigma_1^{-2}$ & $\sigma_1^{2}\sigma_2^{3}\sigma_3$  & \\
\hline
$\sigma_1^{2}\sigma_2^{2}\sigma_3^{2}$ & $\sigma_1\sigma_2\sigma_1\sigma_3\sigma_2^{2}\sigma_3\sigma_2\sigma_1\sigma_2^{-1}\sigma_1^{-2}$ & $\sigma_1^{3}\sigma_2\sigma_1\sigma_3\sigma_2\sigma_1^{-1}$ & $\sigma_1^{2}\sigma_2^{2}\sigma_3^{2}$   \\
\hline
$\sigma_1^{2}\sigma_2\sigma_3^{2}\sigma_2$ & $\sigma_3$ & $\sigma_1\sigma_2\sigma_1^{-1}$ & $\sigma_1^{2}\sigma_2\sigma_3^{2}\sigma_2$   \\
\hline
$\sigma_1\sigma_2^{4}\sigma_3$ & $\sigma_1\sigma_2^{3}\sigma_3^{-1}\sigma_1^{-1}\sigma_2^{-1}\sigma_1^{-1}$ & $\sigma_1^{2}\sigma_2\sigma_1\sigma_3\sigma_2$ & \\
\hline
\end{tabular}
}

\vspace{.5cm}
 Actually, every time that we tried to reduce the number of generators associated to a conjugacy class,
we were able to keep just $n-1$. Remark also that there are 1634 different conjugacy classes of elements
of length $l$ $(4\leq l \leq 20)$ in $B_3^+$, all of them with no
more than two generators. So all these evidences led us to think that the centralizer of every braid in
$B_n$ could be generated by at most $n-1$ elements.

 As we said, this conjecture turned out to be false, since a family of counterexamples due to N. V. Ivanov
gives a lower bound for the number of generators which is a quadratic in $n$. Precisely, there has been
recently shown~\cite{GW} that the centralizer of every element in $B_n$ can
be generated by less than $\frac{n(n-1)}{2}$ elements.

 In any case, the above results are valid just for braids, so we still would like to have an upper bound for
the minimal number of generators of $Z(a)$, in the general case of Garside groups.

\vspace{.5cm}
\noindent
{\bf Acknowledgements:} The authors want to thank the `Laboratorie de Topologie de l'Universit\'e
de Bourgogne', where we started to work in this subject, and to Luis Paris, Alain Jacquemard,
Jos\'e Mar\'{\i}a Tornero, Carmen Le\'on, Mustafa Korkmaz, Arkadius Kalka and Bert Wiest for their valuable help.

\noindent {\footnotesize Nuno Franco: \\
Dep. de Matem\'atica, CIMA-UE, Universidade de \'Evora,
7000-\'Evora (PORTUGAL), E-mail: {\em nmf@uevora.pt} \\
Universit\'e de Bourgogne, Laboratoire de Topologie, UMR 5584 du CNRS, B.P. 47870,
21078-Dijon Cedex (FRANCE). \vspace{.2cm}\\
Juan Gonz\'alez-Meneses: \\
Dep. Matem\'atica Aplicada I, ETS Arquitectura, Univ. de Sevilla,
Av. Reina Mercedes 2, 41012-Sevilla (SPAIN). E-mail: {\em meneses@us.es}}
\end{document}

%% file: Delta475.pstex_t
\begin{picture}(0,0)%
\epsfig{file=Delta475.pstex}%
\end{picture}%
\setlength{\unitlength}{3108sp}%
\begingroup\makeatletter\ifx\SetFigFont\undefined%
\gdef\SetFigFont#1#2#3#4#5{%
  \reset@font\fontsize{#1}{#2pt}%
  \fontfamily{#3}\fontseries{#4}\fontshape{#5}%
  \selectfont}%
\fi\endgroup%
\begin{picture}(1465,2082)(1281,-1423)
\put(1306,524){\makebox(0,0)[lb]{\smash{\SetFigFont{9}{10.8}{\rmdefault}{\mddefault}{\updefault}\special{ps: gsave 0 0 0 setrgbcolor}1\special{ps: grestore}}}}
\put(1756,524){\makebox(0,0)[lb]{\smash{\SetFigFont{9}{10.8}{\rmdefault}{\mddefault}{\updefault}\special{ps: gsave 0 0 0 setrgbcolor}2\special{ps: grestore}}}}
\put(2206,524){\makebox(0,0)[lb]{\smash{\SetFigFont{9}{10.8}{\rmdefault}{\mddefault}{\updefault}\special{ps: gsave 0 0 0 setrgbcolor}3\special{ps: grestore}}}}
\put(2656,524){\makebox(0,0)[lb]{\smash{\SetFigFont{9}{10.8}{\rmdefault}{\mddefault}{\updefault}\special{ps: gsave 0 0 0 setrgbcolor}4\special{ps: grestore}}}}
\end{picture}

%% file: lattice.pstex_t
\begin{picture}(0,0)%
\epsfig{file=lattice.pstex}%
\end{picture}%
\setlength{\unitlength}{3108sp}%
\begingroup\makeatletter\ifx\SetFigFont\undefined%
\gdef\SetFigFont#1#2#3#4#5{%
  \reset@font\fontsize{#1}{#2pt}%
  \fontfamily{#3}\fontseries{#4}\fontshape{#5}%
  \selectfont}%
\fi\endgroup%
\begin{picture}(7449,4119)(439,-4213)
\put(6751,-1546){\makebox(0,0)[lb]{\smash{\SetFigFont{9}{10.8}{\rmdefault}{\mddefault}{\updefault}$\sigma_2\sigma_3\sigma_2\sigma_1$}}}
\put(5266,-1546){\makebox(0,0)[lb]{\smash{\SetFigFont{9}{10.8}{\rmdefault}{\mddefault}{\updefault}$\sigma_2\sigma_1\sigma_3\sigma_2$}}}
\put(3781,-1546){\makebox(0,0)[lb]{\smash{\SetFigFont{9}{10.8}{\rmdefault}{\mddefault}{\updefault}$\sigma_1\sigma_3\sigma_2\sigma_1$}}}
\put(2296,-1546){\makebox(0,0)[lb]{\smash{\SetFigFont{9}{10.8}{\rmdefault}{\mddefault}{\updefault}$\sigma_1\sigma_2\sigma_3\sigma_2$}}}
\put(856,-1546){\makebox(0,0)[lb]{\smash{\SetFigFont{9}{10.8}{\rmdefault}{\mddefault}{\updefault}$\sigma_1\sigma_2\sigma_1\sigma_3$}}}
\put(3646,-916){\makebox(0,0)[lb]{\smash{\SetFigFont{9}{10.8}{\rmdefault}{\mddefault}{\updefault}$\sigma_1\sigma_2\sigma_3\sigma_2\sigma_1$}}}
\put(5221,-916){\makebox(0,0)[lb]{\smash{\SetFigFont{9}{10.8}{\rmdefault}{\mddefault}{\updefault}$\sigma_2\sigma_1\sigma_3\sigma_2\sigma_1$}}}
\put(7111,-2176){\makebox(0,0)[lb]{\smash{\SetFigFont{9}{10.8}{\rmdefault}{\mddefault}{\updefault}$\sigma_3\sigma_2\sigma_1$}}}
\put(4051,-331){\makebox(0,0)[lb]{\smash{\SetFigFont{12}{14.4}{\rmdefault}{\mddefault}{\updefault}$\Delta$}}}
\put(5806,-2176){\makebox(0,0)[lb]{\smash{\SetFigFont{9}{10.8}{\rmdefault}{\mddefault}{\updefault}$\sigma_2\sigma_3\sigma_2$}}}
\put(4546,-2176){\makebox(0,0)[lb]{\smash{\SetFigFont{9}{10.8}{\rmdefault}{\mddefault}{\updefault}$\sigma_2\sigma_1\sigma_3$}}}
\put(3286,-2176){\makebox(0,0)[lb]{\smash{\SetFigFont{9}{10.8}{\rmdefault}{\mddefault}{\updefault}$\sigma_1\sigma_3\sigma_2$}}}
\put(2026,-2176){\makebox(0,0)[lb]{\smash{\SetFigFont{9}{10.8}{\rmdefault}{\mddefault}{\updefault}$\sigma_1\sigma_2\sigma_3$}}}
\put(676,-2176){\makebox(0,0)[lb]{\smash{\SetFigFont{9}{10.8}{\rmdefault}{\mddefault}{\updefault}$\sigma_1\sigma_2\sigma_1$}}}
\put(6931,-2806){\makebox(0,0)[lb]{\smash{\SetFigFont{9}{10.8}{\rmdefault}{\mddefault}{\updefault}$\sigma_3\sigma_2$}}}
\put(5446,-2806){\makebox(0,0)[lb]{\smash{\SetFigFont{9}{10.8}{\rmdefault}{\mddefault}{\updefault}$\sigma_2\sigma_3$}}}
\put(3961,-2806){\makebox(0,0)[lb]{\smash{\SetFigFont{9}{10.8}{\rmdefault}{\mddefault}{\updefault}$\sigma_2\sigma_1$}}}
\put(2476,-2806){\makebox(0,0)[lb]{\smash{\SetFigFont{9}{10.8}{\rmdefault}{\mddefault}{\updefault}$\sigma_1\sigma_3$}}}
\put(1036,-2806){\makebox(0,0)[lb]{\smash{\SetFigFont{9}{10.8}{\rmdefault}{\mddefault}{\updefault}$\sigma_1\sigma_2$}}}
\put(5536,-3436){\makebox(0,0)[lb]{\smash{\SetFigFont{9}{10.8}{\rmdefault}{\mddefault}{\updefault}$\sigma_3$}}}
\put(4051,-3436){\makebox(0,0)[lb]{\smash{\SetFigFont{9}{10.8}{\rmdefault}{\mddefault}{\updefault}$\sigma_2$}}}
\put(2611,-3436){\makebox(0,0)[lb]{\smash{\SetFigFont{9}{10.8}{\rmdefault}{\mddefault}{\updefault}$\sigma_1$}}}
\put(4096,-4111){\makebox(0,0)[lb]{\smash{\SetFigFont{12}{14.4}{\rmdefault}{\mddefault}{\updefault}1}}}
\put(2116,-916){\makebox(0,0)[lb]{\smash{\SetFigFont{9}{10.8}{\rmdefault}{\mddefault}{\updefault}$\sigma_1\sigma_2\sigma_1\sigma_3\sigma_2$}}}
\end{picture}

%% file: MCG1.pstex_t
\begin{picture}(0,0)%
\epsfig{file=MCG1.pstex}%
\end{picture}%
\setlength{\unitlength}{4144sp}%
\begingroup\makeatletter\ifx\SetFigFont\undefined%
\gdef\SetFigFont#1#2#3#4#5{%
  \reset@font\fontsize{#1}{#2pt}%
  \fontfamily{#3}\fontseries{#4}\fontshape{#5}%
  \selectfont}%
\fi\endgroup%
\begin{picture}(4725,1171)(271,-920)
\put(1666,-241){\makebox(0,0)[lb]{\smash{\SetFigFont{10}{12.0}{\rmdefault}{\mddefault}{\updefault}\special{ps: gsave 0 0 0 setrgbcolor}$\sigma_2\sigma_1$\special{ps: grestore}}}}
\put(1666,-781){\makebox(0,0)[lb]{\smash{\SetFigFont{10}{12.0}{\rmdefault}{\mddefault}{\updefault}\special{ps: gsave 0 0 0 setrgbcolor}$\sigma_1\sigma_2$\special{ps: grestore}}}}
\put(3376,-241){\makebox(0,0)[lb]{\smash{\SetFigFont{10}{12.0}{\rmdefault}{\mddefault}{\updefault}\special{ps: gsave 0 0 0 setrgbcolor}$\sigma_3\sigma_2$\special{ps: grestore}}}}
\put(3376,-781){\makebox(0,0)[lb]{\smash{\SetFigFont{10}{12.0}{\rmdefault}{\mddefault}{\updefault}\special{ps: gsave 0 0 0 setrgbcolor}$\sigma_2\sigma_3$\special{ps: grestore}}}}
\put(811,-556){\makebox(0,0)[lb]{\smash{\SetFigFont{14}{16.8}{\rmdefault}{\mddefault}{\updefault}\special{ps: gsave 0 0 0 setrgbcolor}$\sigma_1$\special{ps: grestore}}}}
\put(2611,-556){\makebox(0,0)[lb]{\smash{\SetFigFont{14}{16.8}{\rmdefault}{\mddefault}{\updefault}\special{ps: gsave 0 0 0 setrgbcolor}$\sigma_2$\special{ps: grestore}}}}
\put(4411,-556){\makebox(0,0)[lb]{\smash{\SetFigFont{14}{16.8}{\rmdefault}{\mddefault}{\updefault}\special{ps: gsave 0 0 0 setrgbcolor}$\sigma_3$\special{ps: grestore}}}}
\put(4996,-196){\makebox(0,0)[lb]{\smash{\SetFigFont{10}{12.0}{\rmdefault}{\mddefault}{\updefault}\special{ps: gsave 0 0 0 setrgbcolor}$\sigma_1$\special{ps: grestore}}}}
\put(4996,-871){\makebox(0,0)[lb]{\smash{\SetFigFont{10}{12.0}{\rmdefault}{\mddefault}{\updefault}\special{ps: gsave 0 0 0 setrgbcolor}$\sigma_3$\special{ps: grestore}}}}
\put(271,-196){\makebox(0,0)[lb]{\smash{\SetFigFont{10}{12.0}{\rmdefault}{\mddefault}{\updefault}\special{ps: gsave 0 0 0 setrgbcolor}$\sigma_1$\special{ps: grestore}}}}
\put(271,-871){\makebox(0,0)[lb]{\smash{\SetFigFont{10}{12.0}{\rmdefault}{\mddefault}{\updefault}\special{ps: gsave 0 0 0 setrgbcolor}$\sigma_3$\special{ps: grestore}}}}
\put(2611,119){\makebox(0,0)[lb]{\smash{\SetFigFont{10}{12.0}{\rmdefault}{\mddefault}{\updefault}\special{ps: gsave 0 0 0 setrgbcolor}$\sigma_2$\special{ps: grestore}}}}
\end{picture}

%% file: MCG12.pstex_t
\begin{picture}(0,0)%
\epsfig{file=MCG12.pstex}%
\end{picture}%
\setlength{\unitlength}{4144sp}%
\begingroup\makeatletter\ifx\SetFigFont\undefined%
\gdef\SetFigFont#1#2#3#4#5{%
  \reset@font\fontsize{#1}{#2pt}%
  \fontfamily{#3}\fontseries{#4}\fontshape{#5}%
  \selectfont}%
\fi\endgroup%
\begin{picture}(2753,1711)(901,-1775)
\put(2386,-1726){\makebox(0,0)[lb]{\smash{\SetFigFont{10}{12.0}{\rmdefault}{\mddefault}{\updefault}\special{ps: gsave 0 0 0 setrgbcolor}$\sigma_3$\special{ps: grestore}}}}
\put(901,-1006){\makebox(0,0)[lb]{\smash{\SetFigFont{10}{12.0}{\rmdefault}{\mddefault}{\updefault}\special{ps: gsave 0 0 0 setrgbcolor}$\sigma_1\sigma_2\sigma_3$\special{ps: grestore}}}}
\put(2791,-1006){\makebox(0,0)[lb]{\smash{\SetFigFont{10}{12.0}{\rmdefault}{\mddefault}{\updefault}\special{ps: gsave 0 0 0 setrgbcolor}$\sigma_1\sigma_2\sigma_3$\special{ps: grestore}}}}
\put(1351,-1591){\makebox(0,0)[lb]{\smash{\SetFigFont{14}{16.8}{\rmdefault}{\mddefault}{\updefault}\special{ps: gsave 0 0 0 setrgbcolor}$\sigma_2\sigma_3$\special{ps: grestore}}}}
\put(1351,-331){\makebox(0,0)[lb]{\smash{\SetFigFont{14}{16.8}{\rmdefault}{\mddefault}{\updefault}\special{ps: gsave 0 0 0 setrgbcolor}$\sigma_1\sigma_2$\special{ps: grestore}}}}
\put(3196,-331){\makebox(0,0)[lb]{\smash{\SetFigFont{14}{16.8}{\rmdefault}{\mddefault}{\updefault}\special{ps: gsave 0 0 0 setrgbcolor}$\sigma_2\sigma_1$\special{ps: grestore}}}}
\put(3196,-1591){\makebox(0,0)[lb]{\smash{\SetFigFont{14}{16.8}{\rmdefault}{\mddefault}{\updefault}\special{ps: gsave 0 0 0 setrgbcolor}$\sigma_3\sigma_2$\special{ps: grestore}}}}
\put(2386,-196){\makebox(0,0)[lb]{\smash{\SetFigFont{10}{12.0}{\rmdefault}{\mddefault}{\updefault}\special{ps: gsave 0 0 0 setrgbcolor}$\sigma_1$\special{ps: grestore}}}}
\put(2386,-466){\makebox(0,0)[lb]{\smash{\SetFigFont{10}{12.0}{\rmdefault}{\mddefault}{\updefault}\special{ps: gsave 0 0 0 setrgbcolor}$\sigma_2$\special{ps: grestore}}}}
\put(2386,-1456){\makebox(0,0)[lb]{\smash{\SetFigFont{10}{12.0}{\rmdefault}{\mddefault}{\updefault}\special{ps: gsave 0 0 0 setrgbcolor}$\sigma_2$\special{ps: grestore}}}}
\put(3556,-1006){\makebox(0,0)[lb]{\smash{\SetFigFont{10}{12.0}{\rmdefault}{\mddefault}{\updefault}\special{ps: gsave 0 0 0 setrgbcolor}$\sigma_3\sigma_2\sigma_1$\special{ps: grestore}}}}
\put(1711,-1006){\makebox(0,0)[lb]{\smash{\SetFigFont{10}{12.0}{\rmdefault}{\mddefault}{\updefault}\special{ps: gsave 0 0 0 setrgbcolor}$\sigma_3\sigma_2\sigma_1$\special{ps: grestore}}}}
\end{picture}